\documentclass{amsart}
\usepackage{color}
\usepackage[all]{xy}
\newtheorem{example}{Example}
\newtheorem{lem}{Lemma}
\newtheorem*{rmk}{Remark}
\newtheorem{cor}{Corollary}
\newtheorem{thm}{Theorem}

\newcommand{\N}{\mathbb{N}}
\newcommand{\Q}{\mathbb{Q}}
\newcommand{\R}{\mathbb{R}}
\newcommand{\Z}{\mathbb{Z}}
\newcommand{\A}{\alpha}

\newcommand{\E}{\epsilon}

\newcommand{\ST}{\textrm{ : }}

\begin{document}
\title{Heaviness - An Extension of a Lemma of Y. Peres}
\author{David Ralston}
\address{David Ralston \\OSU Department of Mathematics \\ 231 West 18th Ave. \\ Columbus, OH 43210}
\email{ralston@math.ohio-state.edu}
\maketitle

\begin{abstract}
We provide an elementary proof of Y. Peres' lemma on the existence in certain dynamical systems of what we term \emph{heavy points}, points whose ergodic averages consistently dominate the expected value of the ergodic averages.  We also derive several generalizations of Peres' lemma by employing techniques from the simplified proof.
\end{abstract}
\section{The Lemma of Peres and its Immediate Generalizations}
The following lemma is derived by Yuval Peres \cite{peres} using the maximal ergodic theorem.  While this original proof is quite short and natural, we will see that elementary methods yield more general results.  Before proceeding to these broader results, however, the original lemma deserves mention.
\begin{lem}[Peres]{Let $T:X \rightarrow X$ be a continuous transformation of a compact space, and let $\mu$ be a probability measure preserved by $T$.  For every continuous $f:X \rightarrow \R$ there exists some $x \in X$ such that \begin{equation}\label{eqn} \forall N \in \N \hspace{.2 in}\frac{1}{N}\sum_{i=0}^{N-1}f(T^i x) \geq \int_X f d\mu.\end{equation} }\label{peres lemma}
\begin{proof}[Proof (due to Peres)]
  For any $\E>0$, let \[E_{\E}=\left\{ x \in X \ST \forall N \in \N \hspace{.2 in} \frac{1}{N}\sum_{i=0}^{N-1}f(T^ix)\geq \int_X f d\mu -\E \right\}\] and define $g, M$ by \[g=\int_X f d\mu - f -\E , \hspace{.2 in} M(x)=\sup_{N\in \N}\left\{ \frac{1}{N}\sum_{i=0}^{N-1}g(T^i x) \right\},\] so that $E_{\E}=\{x \ST M(x) \leq 0\}$.  The maximal ergodic theorem states that \[0 \leq \int_{M(x)>0}g d\mu = \int_{X \setminus E_{\E}}g d\mu.\] Now, as $\int_X g d\mu = -\E$, we must have that $E_{\E} \neq \emptyset$.  Furthermore, as $f$ is continuous, $E_{\E}$ is closed.  By intersecting the (nested) sets $E_{\E}$ as $\E \rightarrow 0$ in our compact space $X$, we obtain our desired point.
\end{proof}
\end{lem}

A point $x$ which satisfies inequality \eqref{eqn} we term \textit{heavy for} $f$ (under $T$).  If $f=\chi_A$, the characteristic function of a set $A$, we will simply call the point $x$ \textit{heavy for} $A$.  We may also define heaviness in measure spaces without the luxury of a measure-preserving $T$: given any sequence $\{x_i\}_{i=0,1,\ldots}$ in a probability space $\{X,\mu\}$, and a function $f \in L^1(X,\mu)$, we say that the sequence is \textit{heavy for $f$} if the following inequality, similar to \eqref{eqn}, holds:
\begin{equation}\label{hvyseq}\forall N \in \N \hspace{.2 in}\sum_{i=0}^{N-1}f(x_i) \geq N\int_X f d\mu.\end{equation}
In Section \ref{2x 3x}, we will use tools of ergodic theory to investigate one interesting such situation where heaviness is absent, and Example \ref{polynomial} similarly uses a construction of ergodic theory to produce heaviness in another non-obvious setting.
\begin{rmk}
We note that in a nonperiodic ergodic system, if $\mu(A)\notin \Q$, then the set of points heavy for $A$ must be a null set.  Theorem 4 in \cite{halasz} states that the quantity \begin{equation}\label{hal} \frac{1}{N}\sum_{i=0}^{N-1}f(T^i x)-\int_X f d\mu\end{equation} cannot be ultimately strictly positive or ultimately strictly negative except on a null set (provided $T$ is ergodic and nonperiodic), and we omit the proof. \end{rmk}

This remark should not be overextended, however:

\begin{example}{There are nonperiodic ergodic systems and sets $A$ of nontrivial measure ($\mu(A) \notin \{0,1\}$) for which the set of heavy points is of positive measure.}
\begin{proof}  Let $X$ be the space generated by the shift map on the \emph{Morse Sequence} (see \cite{allouche-shallit} for background and further details): \[x=x_1 x_2 x_3 \ldots=01101001\ldots\] The sequence $x$ is formed by the substitution $0 \rightarrow 01$, $1 \rightarrow 10$: \[0 \rightarrow 01 \rightarrow 0110 \rightarrow 01101001 \rightarrow \ldots\]  It is a short exercise to verify that $T$, the shift map $\left(T(x)_n=x_{n+1}\right)$, is uniquely ergodic in the resulting nonatomic system $\left\{ \overline{O^+(x)},T\right\}$, and that $\mu (A)=1/2$, where $A= \left\{ y \ST y_1=1\right\}$ and $\mu$ is the unique preserved measure.  It is quick to show that any word in $B=\{y \ST y_1 y_2 = 11\}$ will be heavy for $A$, and that $\mu (B) >0$.  This system does not contradict the results in \cite{halasz}; in this example, the quantities \eqref{hal} are not ultimately positive on $B$ - they are zero for infinitely many $N$.\end{proof}
\end{example}

\begin{rmk}
  We do not need in Lemma \ref{peres lemma} that $f$ be continuous - upper semi-continuity would imply that the sets $E_{\E}$ are closed.
\end{rmk}

\begin{cor}{Let $X$ be compact and $T:X \rightarrow X$ be continuous.  Let $A \subset X$ be closed.  Then given any $T$-invariant measure $\mu$, there exists some $x \in X$ which is heavy for $A$.}\label{heavyinterval}
\begin{proof}
  The function $\chi_A$ is upper semi-continuous, so we apply Lemma \ref{peres lemma}.
\end{proof}
\end{cor}
Corollary \ref{heavyinterval} seems trivial, but implies the following example:
\begin{example}{Fix $\alpha \in \R$ and $k \in \N$.  Given $\mathbb{T}^k \ni \vec{a}=a_0, a_1, \ldots, a_{k-1}$, let $p^{\vec{a}}(x)=\alpha x^k + a_{k-1}x^{k-1}+\ldots +a_1 x + a_0$.  Let $I\subset S^1$ be a closed set (for example, a connected closed interval).  Then there exists some choice of $\vec{a}$ (a choice of coefficients other than $\A$) such that the sequence $\{p^{\vec{a}}(n)\}_{n=0}^{\infty}$ is heavy for $I$ (satisfies \eqref{hvyseq} with $f =\chi_I$).  Specifically, for all $N$, we have \[\sum_{i=0}^{N-1}\chi_I(p(i)) \geq N\mu (I).\]}\label{polynomial}

\begin{proof}
  The map $T: \mathbb{T}^k \rightarrow \mathbb{T}^k$ given by \[T(x_1, \ldots, x_k)=(x_1+\alpha, x_2+x_1, \ldots, x_k+x_{k-1})\] clearly preserves Lebesgue measure.  Define $p_k^{\vec{a}}(x)=\alpha x^k + a_{k-1}x^{k-1}+\ldots+a_1 x + a_1$ and $p_{i-1}^{\vec{a}}(x)=p_i^{\vec{a}}(x+1)=p_i^{\vec{a}}(x)$.  Note that for any $\vec{x}\in \mathbb{T}^k$, we may find $\vec{a}$ such that $\vec{x}=(p_1^{\vec{a}}(0),\ldots,p_k^{\vec{a}}(0))$.  Then we have (all entries modulo one): \[T^n(p_1^{\vec{a}}(0),\ldots,p_k^{\vec{a}}(0))=(p_1^{\vec{a}}(n),\ldots,p_k^{\vec{a}}(n))\]
  Let $I'=S^1 \times \dots \times S^1 \times I$ and, as $T$ is Lebesgue measure-preserving, and $\chi_{I'}$ is upper semi-continuous (and integrates to $\mu(I)$), there is some point $\vec{x}$ which is heavy for $I'$, which corresponds to a choice of $\vec{a}$ such that the sequence $p_k^{\vec{a}}(n)$ is heavy for $I$.
\end{proof}
\end{example}
\begin{rmk}
It is possible to derive this result (with a few arguments) from knowing the equidistribution of $p^{\vec{a}}(n)$ when $\alpha \notin \Q$, but note that our proof relied only on our map $T$ being measure-preserving, and was therefore independent of the irrationality of $\alpha$, while equidistribution requires knowledge of $T$ being uniquely ergodic, a stronger result which requires $\alpha \notin \Q$.  For a proof of unique ergodicity of $T$, and therefore equidistribution of $p^{\vec{a}}(n)$ for all choices of $\vec{a}$ (when $\alpha \notin \Q$), we refer the reader to \cite{furstenberg81}.
\end{rmk}

In the next section we formalize our terminology and notation, and we obtain a purely measure-theoretical version of Peres' Lemma (Lemma \ref{peres lemma}).  If we restrict our attention to finding only points whose ergodic averages dominate the expected value through some finite time $N$, we may disregard all topological information, and even finiteness of $\mu$ in specific cases, and find a nontrivial set of such points.  The third section is devoted to extending these results for invertible actions to negative times.  Finally, the fourth section addresses further connections to the realm of heavy sequences, using dynamical systems techniques to illustrate an interesting counterexample.

\section{Definitions and the First Theorem}

Now let $\{X, \mu, T\}$ be a probability-measure-preserving system.  Let $f \in L^1(X,\mu)$.  The \emph{$N^{\textit{th}}$ partial sum} of a point $x$ is given by $S_N(x)=\sum_{i=0}^{N-1}f(T^ix)$.  Define $\psi(x)$ by \[\psi(x)=\inf_{n \in N}\left\{n \ST S_n(x)-n\int_Xf d\mu <0\right\}\] and then partition the space $X$ into the countable collection of sets $\mathcal{E}(n)=\{x \ST \psi(x)= n\}$.  If $\psi(x) > n$, we say that $x$ is \emph{heavy through time $n$}, and we denote such points by \[\mathcal{H}(n)=\bigcup_{i=n+1}^{\infty}\mathcal{E}(i)\] (where the above union should be understood to include the set $\mathcal{E}(\infty)$) and if $\psi(x)=\infty$, we simply say that $x$ is \emph{heavy}.  We denote (here there is no $\mathcal{H}(\infty)$ to include in the intersection): \[\mathcal{H}=\mathcal{E}(\infty)=\bigcap_{i=1}^{\infty}\mathcal{H}(i).\]

The relation to the lemma of Peres is immediate: the sets $\mathcal{H}(n)$ represent points whose ergodic averages through time $n$ dominate the expected value.  Then Lemma \ref{peres lemma} becomes:
\begin{lem}{If $T:X \rightarrow X$ is a continuous transformation on a compact space preserving a probability measure $\mu$, and $f \in L^1(X,\mu)$ is upper semi-continuous, then $\mathcal{H} \neq \emptyset$.}
\end{lem}

By approaching Lemma \ref{peres lemma} from a finite perspective (the sets $\mathcal{H}(n)$), however, it is possible to gain more information, and without needing topological data or the Maximal Ergodic Theorem.  Some elements of the following proof are similar to the proof of the Maximal Ergodic Theorem due to Y. Katznelson and D. Ornstein, given by K. Petersen in his book, \cite[pp 27-30]{petersen}.  We now assume that $\int_X f d\mu=0$: if $\mu(X)=1$, then we may adjust $f$ by a constant to meet this criteria.  However, by forcing $\int_X f d\mu =0$, we are able to drop the condition that $\mu(X)<\infty$, and still derive the following:

\begin{thm}{Let $\{X,\mu \}$ be a measure space, $T$ an action which preserves $\mu$, and
let $f \in L^1(X,\mu)$ be such that $\int_X f(x)d\mu=0$.  Then
$\forall N \geq 0$, $\mu \left(\mathcal{H}(N)\right)>0$.  That is, the set of points which are heavy through time $N$
is of positive measure.}\label{generalexistence}
\begin{proof}
Without loss of generality, let $T$ be invertible (otherwise, replace our system with its natural invertible extension, as outlined in the book of Cornfeld, Fomin, and Sinai \cite[pp 239-241]{cornfeld-sinai}).  Assume that we may partition our space $X$ up to measure zero into sets $\mathcal{E}(1), \ldots, \mathcal{E}(N)$.  Note that $x\in \mathcal{E}(i) \Rightarrow S_{i}(x)<0$.  Set $n_1=N$, and let $A_{1,0}=\mathcal{E}(N)$ and $T^i A_{1,0}=A_{1,i}$ for $i=0,1,\ldots,n_1-1$.  Now, let $X_2=X\setminus B_1$, where $B_1=\bigcup_{i=0}^n A_{1,i}$.  Let \[n_2=\max \{i\in \N \ST \mu\left(\mathcal{E}(i)\cap X_2\right)>0 \},\]  \[A_{2,0}=\mathcal{E}(n_2)\cap X_2, \hspace{.1 in}A_{2,i}=T^i A_{2,0} \textrm{ for }i=1,2,\ldots n_2-1.\]  Similarly define $A_{3,i}$, $A_{4,i}$, etc.  As $n_{i+1}<n_i$, this process terminates at some stage $m$.  See Figure \ref{generalexistencepic}.

As outlined in Petersen \cite[p 29]{petersen}, the sets $A_{i,j}$ are all disjoint.  This fact may be deduced by first observing that $A_{1,0}\cap A_{1,i}$ must be empty (a point $x$ in the intersection would have to be heavy through time $n_1+i-1>n_1-1$, contradicting $x \in A_{1,0}$), and applying similar reasoning to other possibilities.  Furthermore:
\[\int_{\bigcup_{i=0}^{n_j-1} A_{j,i}} f(x) d\mu = \sum_{i=0}^{n_j-1}\int_{A_{j,i}}f(x) d\mu=\int_{A_{j,0}}\sum_{i=0}^{n_j-1}f(T^i x)d\mu=\int_{A_{j,0}}S_{n_j}(x)d\mu<0\]

We have created a partition of $X$ (up to a null set) into disjoint sets of nonzero measure (the rows $\bigcup_{i=0}^{n_k-1}A(k,i)$ in Figure \ref{generalexistencepic}), over each of which the function $f$ integrates to \emph{strictly less} than zero, contradicting the fact that $\int_X f d\mu=0$.
\end{proof}

\begin{figure}[h]
\caption{\label{generalexistencepic} A partition of the space $X$ into disjoint sets, where each row integrates to strictly less than zero. }\center{\makebox{
\xymatrix@C=35pt@R=15pt{
*+[F-,]{A_{1,0}}\ar[r]^{T} & *+[F-,]{A_{1,1}}\ar[r]^{T}  & \dots \ar[r]^T  & *+[F-,]{A_{1,n_1-2}}\ar[r]^T & *+[F-,]{A_{1,n_1-1}} \\
*+[F-,]{A_{2,0}}\ar[r]^{T} & \dots \ar[r]^T       & *+[F-,]{A_{2,n_2-2}} \ar[r]^{T} & *+[F-,]{A_{2,n_2-1}}  &  \\
\vdots      & \vdots  &  \vdots    &  & &\\
*+[F-,]{A_{m,0}}\ar[r]^{T} &\dots \ar[r]^T &  *+[F-,]{A_{m,n_m-1}} & & \\
}}}
\end{figure}
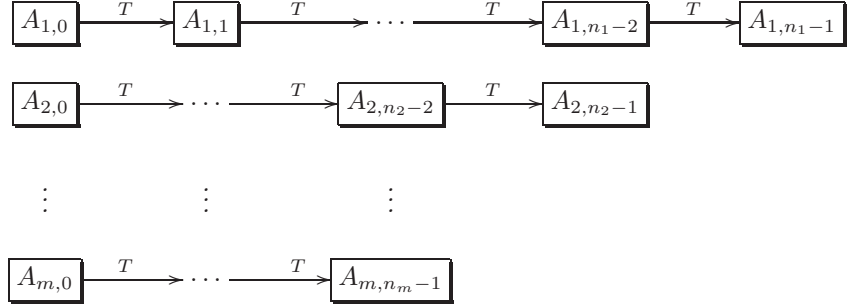

\end{thm}
Note that if $T$ is continuous and our integrable $f$ is upper semi-continuous, then the closed sets $\mathcal{H}(n)$ are nonempty, as they are of positive measure.  Further assuming $X$ to be compact, then, we derive Lemma \ref{peres lemma} as a corollary to Theorem ~\ref{generalexistence}.

\section{Generalizations for Negative Times}
\label{all Z}
Recall that for $n \in \N$ we have defined $S_n(x)=f(x)+f(Tx)+\ldots f(T^{n-1}x)$.  For $n \geq 2$, we clearly have the relation that $S_n(x)-S_{n-1}(x)=f(T^{n-1}x)$.  Assuming $T$ is invertible, we extend this relation to derive that $S_0(x)=0$, and for $n < 0$: \[S_n(x)=-\left[f(T^{-1}x)+f(T^{-2}x)+\ldots+f(T^{n}x)\right].\]  Define the sets $\mathcal{H}(n,m)=\{x \ST S_i(x) \geq 0 \textrm{, }i=n,n+1,\ldots,m\},$ and $\mathcal{H}(\Z)=\cap_{n=1}^{\infty}\mathcal{H}(-n,n).$  Note that $\mathcal{H}(0)=X$.  Furthermore, there is the obvious relation that $S_n(T^m x)=S_{n+m}(x)-S_m (x)$.  We immediately conclude the following:
\begin{lem}{Let $m<n$ be integers.  Then $\left( T^m x \in \mathcal{H}(0,n) \right) \Leftrightarrow \left( S_m(x) \leq S_{m+i}(x) \right)$ for $i=0,1,\ldots,n$.  Similarly, $\left( T^m x \in \mathcal{H}(-n,0) \right) \Leftrightarrow \left( S_m(x) \leq S_{m-i}(x)\right)$ for $i=0,1,\ldots n$.}
\end{lem}

By considering $T^{-1}$ and the function $-f$, we reprove Lemma \ref{peres lemma} to claim $\mathcal{H}(-\infty,0)\neq \emptyset$, replacing the assumption of upper semi-continuity with lower semi-continuity.  Similarly, $\mu \left(\mathcal{H}(-n,0)\right) >0$ under the conditions of Theorem \ref{generalexistence}.  Note that if $f$ is continuous, both $\mathcal{H}(-\infty,0)$ and $\mathcal{H}(0,\infty)$ are nonempty.  Is it possible to combine both positive and negative times and retain either Lemma \ref{peres lemma} or Theorem \ref{generalexistence}?  In general, it is not.

\begin{example}{$X=\{0,1\}$, $T=\textrm{Id}$, $f(0)=1$, $f(1)=-1$.}
  Here we have $S_n(0)=n$ for all $n \in \Z$, and $S_n(1)=-n$, so for each point, we may only claim heaviness in one direction (positive or negative times, but not both).
\end{example}

This system is an unnatural example: the space $X$ decomposes into two $T$-invariant closed subspaces.  Under the appropriate conditions, however, we may reach our desired conclusions:

\begin{thm}{Assume $X$ is compact, $f$ is continuous, and $T$ is transitive, invertible, and continuous.  Then $\mathcal{H}(\Z) \neq \emptyset$.  This theorem is an extension of Lemma ~\ref{peres lemma}.}\label{all Z transitive}
\begin{proof}
First, assume we have a point $y\in \mathcal{H}(0,\infty)$ such that $S_n(y)=0$ for infinitely many $n>0$.  Then by taking $x_i=T^{n_i}(y)$ for an increasing sequence of such times, a limit point $x \in \cap_{n=1}^{\infty}\overline{\{x_i\}_{i=n}^{\infty}}$ will work.  We may construct a similar point if there is $z \in \mathcal{H}(-\infty,0)$ with $S_n(z)=0$ for infinitely many $n<0$.  So, without loss of generality, we assume that there is no point.  Therefore, there exist $y \in \mathcal{H}(0,\infty)$ and $z \in \mathcal{H}(-\infty,0)$ with \emph{no} such times, and these points $y$ and $z$ now have the property that for any $N$, there is a small open neighborhood that remains heavy through time $N$ (time $-N$ for the point $z$).

  Let $x_0$ be a point with dense orbit.  Let $U_N \ni y$ be an open neighborhood of $y$ such that $U_N \subset \mathcal{H}(0,N)$, and let $V_N \ni z$ be an open neighborhood of $z$ such that $V_N \subset \mathcal{H}(-N,0)$.  As the orbit of $x_0$ is dense, for every $N>0$, there is an $i_N>0$ and $j_N>i_N+N$ such that $T^{i_N}x \in V_N$ and $T^{j_N}x \in U_N$.  Let $x_N$ be the point $T^{k+1} x_0$, where $i_N\leq k \leq j_N $ is the time where $S_k(x_0)$ is minimized.  Then we find that $x_N \in \mathcal{H}(0,N+j_N-k) \subset \mathcal{H}(0,N)$, and also $x_N \in \mathcal{H}(-N-k+i_N,0) \subset \mathcal{H}(-N,0)$.  By letting $x \in \cap_{n=1}^{\infty}\overline{\{x_i \}_{i=n}^{\infty}}$, we find that this $x$ has the desired property.
\end{proof}

\end{thm}

 Note that all we needed to complete the proof was a point whose orbit was guaranteed to land in the proper $\mathcal{H}(0,N)$ and $\mathcal{H}(-N,0)$.  Under the assumptions of Theorem \ref{generalexistence}, both of these sets are of positive measure for every $N$.  If we assume ergodicity of $T$, then almost every orbit has this property, and therefore we conclude the following:

\begin{thm}{Let $T$ be invertible and preserving probability measure $\mu$.  Then $T$ is ergodic if and only if for every integrable $f$ and $n_1<n_2$, the set $\mathcal{H}(n_1,n_2)$ is of positive measure.  This theorem is an extension of Theorem \ref{generalexistence}.}
\begin{proof}
  Assume first that $T$ is ergodic.  If $n_1$ and $n_2$ have the same sign, there is nothing to do but apply Theorem \ref{generalexistence}.  So, assume $n_1<0<n_2$.  Almost every point orbits into this set (almost every point orbits into $\mathcal{H}(n_1,0)$ and later into $\mathcal{H}(0,n_2)$, so somewhere in between it must have orbited into $\mathcal{H}(n_1,n_2)$, analogously to Lemma \ref{all Z transitive}), and therefore the set cannot be of zero measure.

  Ergodicity is very important, however, in assuring that almost every point orbits into the desired sets.  Let $A$ be such that $1>\mu A>0$ and $A$ is $T$-invariant.  Then letting $f(x)=\chi_A(x)-\mu A$, the only points heavy for positive times are those points in $A$, and the only points heavy for negative times are those points in $X \setminus A$, so we cannot have $\mu \left(\mathcal{H}(n_1,n_2)\right)>0$ when $n_1<0<n_2$.
\end{proof}
\end{thm}
\section{The sequence $x, 2x, 3x, \ldots$}\label{2x 3x}

We saw in Example \ref{polynomial} that given any closed set $I \subset S^1$ and $\A$, there exists some choice of $\beta$ such that the sequence $n\A +\beta$ is heavy for $I$.  In fact, this $\beta$ is just some element of $\mathcal{H}_{R_{\A}}^{\chi_I}$ (here $R_{\A}$ is the transformation of $S^1=[0,1)$ given by $R_{\A}(x)=x+\A \mod{1}$).  Suppose, however, that we wish to fix $\beta=0$ and allow $\alpha$ to vary?  That is, is there always some point $x$ such that the sequence $x, 2x, 3x, \ldots$ (taken modulo one) is heavy for some closed set $I$?  The answer is `no,' and the proof follows from an application of a $\Z^2$ version of the Rohlin-Halmos `stacking lemma:'

\definecolor{gray}{rgb}{0.9,0.9,0.9}
\begin{figure}[hbt]
\caption{\label{bosh example} Rohlin-Halmos for our two maps $T(x)=2x \mod{1}$ and $S(x)=3x \mod{1}$: each set is disjoint from the others, and altogether they fill the space up to $\E$.  The selected set $B$ (a `checkerboard' subset of the $N \times N$ subsquare) will be disjoint from $T(B)$ and $S(B)$.  In this example, $N$ is odd.}\centering
\makebox{
\xymatrix{\\
\fbox{$t^N A$}\ar[d]^{T} & \fcolorbox{black}{gray}{$t^N s A$}\ar[d]^{T}\ar[l]^{S}  & \dots \ar[d]^{T}\ar[l]^{S}  & \fbox{$t^N s^{N-1} A$}\ar[d]^{T}\ar[l]^{S} & \fcolorbox{black}{gray}{$t^N s^N A$}\ar[d]^{T}\ar[l]^{S}\\
\fbox{$t^{N-1} A$}\ar[d]^{T} & \fbox{$t^{N-1} s A$}\ar[l]^{S}\ar[d]^{T} & \dots \ar[l]^(.4){S}\ar[d]^{T}  & \fcolorbox{black}{gray}{$t^{N-1} s^{N-1} A$}\ar[l]^(.65){S}\ar[d]^{T} & \fbox{$t^{N-1} s^N A$}\ar[l]^{S}\ar[d]^{T}\\
\vdots\ar[d]^{T} & \vdots\ar[l]^{S}\ar[d]^{T}  & \dots \ar[l]^{S}\ar[d]^{T} & \vdots \ar[l]^{S}\ar[d]^{T} & \vdots\ar[l]^{S}\ar[d]^{T}\\
\fbox{$t A$}\ar[d]^{T} & \fcolorbox{black}{gray}{$t s A$}\ar[l]^{S}\ar[d]^{T}  & \dots \ar[l]^{S}\ar[d]^{T} & \fbox{$t s^{N-1} A$}\ar[l]^{S}\ar[d]^{T} & \fcolorbox{black}{gray}{$t s^N A$}\ar[l]^{S}\ar[d]^{T}\\
\fbox{$A$} & \fbox{$s A$}\ar[l]^{S} & \dots \ar[l]^{S}  & \fbox{$s^{N-1} A$}\ar[l]^{S} & \fbox{$s^N A$}\ar[l]^{S}\\
}}
\end{figure}
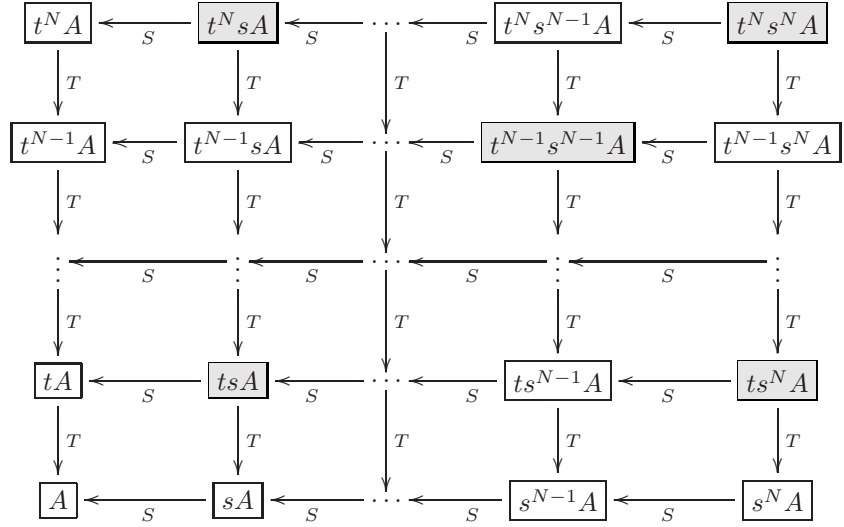

\begin{example}{There is some $I$, a finite collection of closed intervals, such that no $x \in S^1$ such that $x, 2x, 3x, \ldots$ is heavy for I.}
\begin{proof}
  The proof rests on a $\Z^d$ version of the classical Rohlin-Halmos ``stacking lemma," given by Y. Katznelson and B. Weiss \cite{katznelson-weiss}.  Let $T(x)=2x \mod{1}$, and $S(x)=3x \mod{1}$.  These two maps clearly commute with one another and preserve Lebesgue measure.  Furthermore, for any $m,n$, the set of $x$ such that $T^n S^m x=x$ is of zero measure (necessarily being rational), so the group action generated in the corresponding natural extension (again, see \cite[pp 239-241]{cornfeld-sinai}) is totally aperiodic.  For convenience, we will write $t=T^{-1}$ and $s=S^{-1}$.

  For any $N$ and $\E >0$, by the $\Z^2$ version of the ``stacking lemma" referenced above, there is some set $A$, without loss of generality a finite union of closed intervals, such that for all $0\leq n,m \leq N$, the sets $t^{n}s^{m}A$ are all pairwise disjoint, and \[\mu \left( X \setminus \left(\bigcup_{n,m=0}^{N} t^n s^m A\right)\right)<\E.\]  Define \[B=\bigcup_{\substack{2|(p+q) \\N \geq p,q >0}}t^p s^q A\] and see that $\mu (B)= C \mu A$, where $C=N^2/2$ or $(N^2+1)/2$, depending on $N=0 $ or $1 \mod{2}$.  Also, $B$ is a finite union of closed intervals.  For $N=5$ and $\E<1/13$, then, $\mu(B)>1/3$.  By construction, $T(B)\cap B=S(B)\cap B = \emptyset$, so for any $x \in B$, $2x$ and $3x$ both belong to $X\setminus A$ (here the $p+q=0\mod{2}$, $p,q>0$ condition is used).  Therefore, for $x \notin B$, the sequence $\{ix\}_{i=1,2,\ldots}$ violates inequality \eqref{hvyseq} for $N=1$, and for $x \in B$, inequality \eqref{hvyseq} fails for $N=3$.

  See Figure \ref{bosh example} for the construction of $B$.  Note that such a set $B$ as given in the diagram is difficult to construct explicitly, consisting of very many intervals of very small size.  Also note that while we used the maps $T$ and $S$, our end result is a statement on the heaviness properties of the non-dynamic sequence $\{nx\}_{n=1}^{\infty}$.
\end{proof}
\end{example}

While this example shows that we may not claim that \emph{any} closed set $A$ has some point $x \in A$ such that $\{ix\}_{i=1,2,\ldots}$ is heavy for $A$, for certain choices of $A$, the set of such $x$ is of positive Hausdorff dimension.  An investigation of this setting is the topic of a forthcoming paper, joint with M. Boshernitzan \cite{bosh-ralston1}.  In particular, we develop the surprising characterization of the set of $x$ such that $\{ix\}_{i=1}^{\infty}$ (taken modulo one) is heavy for the interval $[0,\frac{1}{k})$: it is exactly those $x$ such that all odd-indexed entries of the continued fraction expansion for $x$ are divisible by $k$ (where the continued fraction terminates at an even index, if $x \in \Q$).

\section{Acknowledgements}
The author is greatly indebted to Michael Boshernitzan, who first initiated the study of heaviness in the context of the sequence $x, 2x, 3x, \ldots$ in \cite{bosh-ralston1}, and began transferring the ideas to dynamical systems.  Karl Petersen was very helpful in providing references to related works, notably \cite{peres}.  William Veech provided the reference to \cite{halasz}, and David Damanik, Jon Chaika, and Heather Hardway have all been patient listeners in the development of the ideas present.


\begin{thebibliography}{9}

\bibitem{allouche-shallit}
J. P. Allouche and J. Shallit, \emph{The Ubiquitous Prouhet-Thue-Morse Sequence}, Sequences and Their Applications (Proc. SETA 1998), 1999.

\bibitem{bosh-ralston1}
M. Boshernitzan and D. Ralston, \emph{Continued Fractions and Heavy Sequences}, Proc. Amer. Math. Soc., \textbf{137} (2009), 3177-3185.

\bibitem{cornfeld-sinai}
I.P. Cornfeld, S.V. Fomin and Ya.G. Sinai, \emph{Ergodic Theory}, Springer-Verlag, 1982.

\bibitem{furstenberg81}
H. Furstenberg, \emph{Recurrence in Ergodic Theory and Combinatorial Number Theory}, Princeton University Press, 1981.

\bibitem{halasz}
G. Hal\'{a}sz, \emph{Remarks on the Remainder in Birkhoff's Ergodic Theorem}, Acta Mathematica \textbf{28}, 1976, 389--395.

\bibitem{katznelson-weiss}
Y. Katznelson and B. Weiss, \emph{Commuting Measure-Preserving Transformations}, Israel Journal of Mathematics \textbf{12}, 1972, 161--173.

\bibitem{peres}
Y. Peres, \emph{A Combinatorial Application of the Maximal Ergodic Theorem}, Bulletin of the London Mathematical Society \textbf{20:3}, 1988, 248--252.

\bibitem{petersen}
K. Petersen, \emph{Ergodic Theory}, Cambridge University Press, 1983.

\end{thebibliography}
\end{document}